\documentclass[12pt]{amsart}
\usepackage[T1]{fontenc}              
\usepackage[latin1]{inputenc}         
\usepackage{times}                    

\makeatletter

\newcommand{\LyX}{L\kern-.1667em\lower.25em\hbox{Y}\kern-.125emX\spacefactor1000}

\theoremstyle{plain}    
\newtheorem{thm}{Theorem}[section]
\numberwithin{equation}{section} 
\numberwithin{figure}{section} 
\theoremstyle{plain}    
\newtheorem{cor}[thm]{Corollary} 
\theoremstyle{plain}    
\newtheorem{lem}[thm]{Lemma} 
\theoremstyle{plain}    
\newtheorem{prop}[thm]{Proposition} 
\theoremstyle{remark}    
\newtheorem*{acknowledgement*}{Acknowledgement} 

\def\Tr{\operatorname{Tr}}
\def\III{\operatorname{III}}
\def\II{\operatorname{II}}
\usepackage{pslatex}

\makeatother

\begin{document}

\title{Prime Type III Factors.}

\author{Dimitri Shlyakhtenko}

\address{Department of Mathematics, UCLA, Los Angeles, CA 90095}

\email{shlyakht@math.ucla.edu}

\date{\today}

\thanks{Research supported by an NSF postdoctoral fellowship.}

\begin{abstract}
We show that for each \( 0<\lambda <1 \), the free Araki-Woods factor of type
III\( _{\lambda } \) cannot be written as a tensor product of two diffuse von
Neumann algebras (i.e., is prime), and does not contain a Cartan subalgebra.
\end{abstract}
\maketitle

\section{Introduction.}

A von Neumann algebra \( M \) is called \emph{prime}, if it cannot be written
as a tensor product of two diffuse von Neumann algebras. Using Voiculescu's
free entropy theory \cite{dvv:entropy3}, Ge \cite{ge:entropy2} and later Stefan
\cite{stephan:thinness} gave examples of prime factors of type II\( _{1} \)
(and hence of type II\( _{\infty } \)). We give an example of a separable prime
factor of type III: we show that for each \( 0<\lambda <1 \), the type III\( _{\lambda } \)
free Araki-Woods factor \( T_{\lambda } \) introduced in \cite{shlyakht:quasifree:big}
is prime. The main idea of the proof is to interpret the decomposition \( T_{\lambda }=A\otimes B \)
as a condition on its core. We then use Stefan's result \cite{stephan:thinness}
showing that \( L(\mathbb {F}_{\infty }) \) cannot be written as the closure
of the linear span of \( N\cdot C_{1}\cdot C_{2} \) where \( N \) is a II\( _{1} \)
factor, which is not prime, and \( C_{i} \) are abelian von Neumann algebras. 

We also prove existence of separable type III factors that do not have Cartan
subalgebras by showing that \( T_{\lambda } \), \( 0<\lambda <1 \) has no
Cartan subalgebras. The key ingredient is Voiculescu's result on the absence
of Cartan subalgebras in \( L(\mathbb {F}_{\infty }) \) \cite{dvv:entropy3}.

\begin{acknowledgement*}
This work was carried out while visiting Centre Émile Borel, Institut Henri
Poicaré, Paris, France, to which I am grateful for the friendly and encouraging
atmosphere. I would like to especially thank the organizers of the Free Probability
and Operator Spaces program at IHP, Professors P. Biane, G. Pisier and D. Voiculescu,
for a very stimulating semester. I would also like to thank M. Stefan for many
useful conversations.
\end{acknowledgement*}

\section{\protect\( T_{\lambda }\protect \) is prime.}

Recall first the following theorem, due to Connes (see \cite{connes}, Sections
4.2 and 4.3):

\begin{thm}
\label{thrm:connesresultes}Let \( M \) be a separable type \( \III _{\lambda } \)
factor with \( 0<\lambda <1 \). Then there exists a faithful normal state \( \phi  \)
on \( M \), with the following properties:
\begin{enumerate}
\item The centralizer \( M^{\phi }=\{m\in M:\phi (mn)=\phi (nm)\quad \forall n\in M\} \)
is a factor of type \( \operatorname {II}_{1} \);
\item The modular group \( \sigma ^{\phi }_{t} \) of \( \phi  \) is periodic, of
period exactly \( 2\pi /\log \lambda  \);
\item \( M \) is generated as a von Neumann algebra by \( M^{\phi } \) and an isometry
\( V \), satisfying:

\begin{enumerate}
\item \( V^{*}V=1 \), \( V^{k}(V^{*})^{k}\in M^{\phi } \) for all \( k \);
\item \( \sigma ^{\phi }_{t}(V)=\lambda ^{-it}(V) \); in particular, \( \phi (V^{k}(V^{*})^{k})=\lambda ^{k}\phi ((V^{*})^{k}V^{k})=\lambda ^{k}\phi (1)=\lambda ^{k} \);
\item \( V \) normalizes \( M^{\phi } \): \( VmV^{*} \) and \( V^{*}mV \) are
both in \( M^{\phi } \) if \( m\in M^{\phi } \).
\end{enumerate}
\end{enumerate}
The weight \( \phi \otimes \Tr (B(\ell ^{2})) \) is unique up to scalar multiples
and up to conjugation by (inner) automorphisms of \( M\cong M\otimes B(\ell ^{2}) \).
Moreover, property (2) implies (1) and (3) and property (1) implies (2) and
(3). In particular, if \( \phi _{1} \) and \( \phi _{2} \) satisfy either
(1) or (2), the centralizers \( M^{\phi _{1}} \) and \( M^{\phi _{2}} \) are
stably isomorphic: \( M^{\phi _{1}}\otimes B(\ell ^{2})\cong M^{\phi _{2}}\otimes B(\ell ^{2}) \).
\end{thm}
The existence of such a state can be easily seen by writing \( M \) as the
crossed product of a type II\( _{\infty } \) factor \( C \) by a trace-scaling
action of \( \mathbb {Z} \): set \( \hat{\phi } \) to be the crossed-product
weight (where \( C \) is taken with its semifinite trace). Next, compress to
a finite projection \( p\in C \) and set \( \phi =\hat{\phi }(p\cdot p) \).
The isometry \( V \) is precisely the compression of the unitary \( U \),
implementing the trace-scaling action of \( \mathbb {Z} \).

Recall that a von Neumann algebra \( M \) is called \emph{full}, if its group
of inner automorphisms is closed in the \( u \)-topology inside its group of
all automorphisms (see \cite{connes:full}).

\begin{lem}
\label{lemma:decomposition}Let \( M \) be a full type \( \III _{\lambda } \)
factor. Assume that \( M=A_{1}\otimes A_{2} \), where \( A_{1} \) and \( A_{2} \)
are von Neumann algebras. Then \( A_{1} \) and \( A_{2} \) are full both factors,
and exactly one of the following must hold true:
\begin{enumerate}
\item \( A_{1} \) and \( A_{2} \) are both of type \( \III _{\lambda _{1}} \) and
\( \III _{\lambda _{2}} \), respectively, and \( \lambda _{1},\lambda _{2} \)
satisfy: (i) \( 0<\lambda _{i}<1 \), \( i=1,2 \), (ii) \( \lambda _{1}^{\mathbb {Z}}\cap \lambda _{2}^{\mathbb {Z}}=\lambda ^{\mathbb {Z}} \)
;
\item For some \( i\neq j \), \( A_{i} \) is of type \( \III _{\lambda } \) and
\( A_{j} \) is of type \( \operatorname {II} \);
\item For some \( i\neq j \), \( A_{i} \) is of type \( \III _{\lambda } \) and
\( A_{j} \) is of type \( \operatorname {I} \).
\end{enumerate}
In particular, if we require that \( A_{1} \) and \( A_{2} \) must both be
diffuse, only (1) and (2) can occur. Moreover, if (2) occurs, we may assume
that one of the algebras \( A_{1} \), \( A_{2} \) is of type \( \operatorname {II}_{1} \).
\end{lem}
\begin{proof}
If one of \( A_{1} \), \( A_{2} \) fails to be a factor, then their tensor
product would fail to be a factor, hence both \( A_{1} \) and \( A_{2} \)
must be factors. Similarly, if at least one of \( A_{1} \) and \( A_{2} \)
fails to be full, their tensor product would fail to be full.

If, say, \( A_{1} \) is of type I or type II, then \( A_{2} \) must be type
III, since otherwise \( A_{1}\otimes A_{2} \) would be of type II or type I.
Hence if at least one of \( A_{1} \) and \( A_{2} \) is not type III, the
situation described in (2) or (3) must occur.

If \( A_{1} \) and \( A_{2} \) are both type III, so that \( A \) is type
III\( _{\lambda _{1}} \) and \( A_{2} \) is type III\( _{\lambda _{2}} \),
we must prove that \( \lambda ^{\mathbb {Z}}=\lambda _{1}^{\mathbb {Z}}\cap \lambda _{2}^{\mathbb {Z}} \).
Neither \( \lambda _{1} \) nor \( \lambda _{2} \) can be zero, because then
at least one of \( A_{1} \), \( A_{2} \) would then fail to be full, and hence
\( A_{1}\otimes A_{2} \) would fail to be full. 

Denote by \( T(M) \) the \( T \) invariant of Connes (see \cite{connes},
Section 1.3). Since
\[
\frac{2\pi \mathbb {Z}}{\log \lambda }=T(A_{1}\otimes A_{2})=T(A_{1})\cap T(A_{2})\]
\cite[Theorem  1.3.4(c)]{connes} and \( T(A_{j})=\frac{2\pi \mathbb {Z}}{\log \lambda _{j}} \),
we obtain (1).
\end{proof}
\begin{prop}
\label{label:prop:centralizer}Let \( M \) be a type \( \III _{\lambda } \)
factor, and assume that \( M=A_{1}\otimes A_{2} \), where \( A_{1} \) is a
type \( \III _{\lambda _{1}} \) factor, \( A \) is a type \( \III _{\lambda _{2}} \)
factor, and \( \lambda ^{\mathbb {Z}}=\lambda _{1}^{\mathbb {Z}}\cap \lambda _{2}^{\mathbb {Z}} \).
Let \( \phi _{i} \) be a normal faithful state on \( A_{i} \) as in Theorem
\ref{thrm:connesresultes}, and let \( \phi =\phi _{1}\otimes \phi _{2} \)
be a normal faithful state on \( M \). 

Then the centralizer \( M^{\phi } \) of \( \phi  \) in \( M \) is a factor,
which can be written as a closure of the linear span of \( N\cdot C_{1}\cdot C_{2} \),
where \( N \) is a tensor product of two type \( \II _{1} \) factors, and
\( C_{i} \) are abelian von Neumann algebras. In particular, \( M^{\phi } \)
is not isomorphic to \( L(\mathbb {F}_{\infty }) \).
\end{prop}
\begin{proof}
Since the modular group of \( \phi _{1}\otimes \phi _{2} \) is \( \sigma ^{\phi _{1}}_{t}\otimes \sigma ^{\phi _{2}}_{t} \),
it follows that \( \sigma ^{\phi _{1}\otimes \phi _{2}}_{t} \) has period exactly
\( 2\pi /\log \lambda  \). Hence the centralizer of \( \phi _{1}\otimes \phi _{2} \)
is a factor. 

Choose a decreasing sequence of projections \( p^{(1)}_{k}\in A^{\phi _{1}}_{1} \)
and \( p^{(2)}_{k}\in A^{\phi _{2}}_{2} \), with \( \phi _{i}(p^{(i)}_{k})=\lambda _{i}^{k} \),
and isometries \( V_{i}\in A_{i} \), so that \( V_{i}^{*}V_{i}=1 \), \( V_{i}^{k}(V_{i}^{*})^{k}=p_{k}^{(i)} \),
so that \( V_{i} \) normalizes \( A_{i}^{\phi _{i}} \), and \( A_{i}=W^{*}(A_{i}^{\phi _{i}},V_{i}) \).
Then \( A_{1}\otimes A_{2} \) is densely spanned by elements of the form 
\[
W=V_{1}^{m_{1}}\otimes V_{2}^{n_{1}}\cdot a_{1}^{(1)}\otimes a^{(2)}_{1}\cdot V_{1}^{m_{2}}\otimes V_{2}^{n_{2}}\cdots V_{1}^{m_{k}}\otimes V_{2}^{n_{k}},\quad a_{j}^{(i)}\in A^{\phi _{i}}_{i},\, m_{i},n_{i}\in \mathbb {Z},\]
with the convention that \( V_{i}^{-n}=(V_{i}^{*})^{n} \) if \( n\geq 0 \). 

Using the fact that \( V_{i}^{*}aV_{i},V_{i}aV_{i}^{*}\in A_{i}^{\phi _{i}} \)
whenever \( a\in A_{i}^{\phi _{i}} \), we can rewrite \( W \) as
\[
W=(V^{*}_{1})^{m}\otimes (V_{2}^{*})^{n}\cdot a^{(1)}\otimes a^{(2)}\cdot V_{1}^{l}\otimes V_{2}^{k},\quad a^{(i)}\in A_{i}^{\phi _{i}},\, m,n,k,l\geq 0.\]

Let now \( p_{k}=p_{k}^{(1)}=V_{1}^{k}(V_{1}^{*})^{k}\in A_{1}^{\phi _{1}} \)
be as above. One can choose a diffuse commutative von Neumann algebra \( \mathcal{A} \),
containing \( p_{k} \), \( k\geq 0 \), and so that \( \mathcal{A}\subset A_{1}^{\phi _{1}} \)
and \( V_{1}\mathcal{A}V_{1}^{*},V_{1}^{*}\mathcal{A}V_{1}\subset \mathcal{A} \).
Choose a projection \( \mathcal{A}\ni q_{0}\leq p_{1} \), so that \( q_{0}\perp p_{2} \)
and \( \phi _{1}(q)=N\phi _{1}(1-p_{1})=N(1-\lambda _{1}) \) for some integer
\( N \). Choose projections \( \mathcal{A}\ni q_{1},\dots ,q_{n}\leq 1-p_{1} \),
so that \( \sum _{i=1}^{N}q_{i}=1-p_{1} \) and \( \phi _{1}(q_{i})=\phi _{1}(q_{0})=\frac{1}{N}\phi _{1}(1-p_{1}) \).
Choose matrix units \( \{e_{ij}\}_{0\leq i,j\leq N}\subset A_{1}^{\phi _{1}} \),
so that \( e_{ii}=q_{i} \), \( 0\leq i\leq N \). Let \( C \) be the von Neumann
algebra generated by \( \{V_{1}^{k}e_{ij}(V^{*}_{1})^{k}:1\leq i,j\leq N,k\geq 0\} \).
By our choice of \( e_{ij} \), \( C \) is hyperfinite (notice that \( V^{k}_{1}e_{ii}(V^{*}_{1})^{k}\in \mathcal{A} \),
and \( C \) is in fact the crossed product of \( \mathcal{A}\cong L^{\infty }(X) \)
by a singly-generated equivalence relation). Let \( R_{1}=W^{*}(C,V)\subset A_{1} \).
Then \( R_{1} \) is also hypefinite; in fact, it is the crossed product of
\( C \) by the endomorphism \( x\mapsto V_{1}xV_{1}^{*} \). Notice that \( R_{1} \)
contains \( V_{1} \). Furthermore, for all \( k\geq 0 \), there is a \( d\geq 0 \)
and partial isometries \( r_{1},\dots ,r_{d}\in R_{1}\cap A^{\phi _{1}}_{1} \),
so that \( 1-V_{1}^{k}(V_{1}^{*})^{k}=\sum _{i=1}^{d}r_{i}V^{k}_{1}(V^{*}_{1})^{k}r_{i}^{*} \). 

Construct in a similar way the algebra \( R_{2}\subset A_{2} \), in such a
way that \( V_{2}\in R_{2} \) and for all \( k\geq 0 \), there is a \( d\geq 0 \)
and partial isometries \( r_{1},\dots ,r_{d}\in R_{2}\cap A^{\phi _{2}}_{2} \),
so that \( 1-V_{2}^{k}(V_{2}^{*})^{k}=\sum _{i=1}^{d}r_{i}V^{k}_{2}(V^{*}_{2})^{k}r_{i}^{*} \). 

Notice that \( R_{1}\otimes R_{2}\subset A_{1}\otimes A_{2} \) is globally
fixed by the modular group of \( \phi _{1}\otimes \phi _{2} \). In particular,
this means that \( (R_{1}\otimes R_{2})^{\phi _{1}\otimes \phi _{2}|_{R_{1}\otimes R_{2}}}=(R_{1}\otimes R_{2})\cap (A_{1}\otimes A_{2})^{\phi _{1}\otimes \phi _{2}} \).

Assume now that \( W\in (A_{1}\otimes A_{2})^{\phi _{1}\otimes \phi _{2}} \).
Then \( \sigma _{t}^{\phi _{1}\otimes \phi _{2}}(W)=W \). Hence \( \lambda _{1}^{m-l}\cdot \lambda _{2}^{n-k}=1 \).
It follows that \( W \) can be written in one of the following forms, using
the fact that \( V_{i}^{*}A^{\phi _{i}}_{i}V_{i}\subset A_{i}^{\phi _{i}} \):
either
\[
W=(V_{1}^{*})^{m}\otimes 1\cdot a^{(1)}\otimes a^{(2)}\cdot 1\otimes V_{2}^{k}\]
or
\[
W=1\otimes (V_{2}^{*})^{n}\cdot a^{(1)}\otimes a^{(2)}\cdot V^{l}_{1}\otimes 1,\]
where \( a^{(1)}\in A^{\phi _{1}}_{1} \), \( a^{(2)}\in A_{2}^{\phi _{2}} \)
and \( \lambda _{1}^{m}=\lambda _{2}^{k} \), \( \lambda _{2}^{n}=\lambda _{1}^{l} \).
In the first case, choose \( r_{1},\dots ,r_{d}\in R_{1}\cap A^{\phi _{1}}_{1} \)
for which \( 1-V_{1}^{m}(V_{1}^{*})^{m}=\sum _{i=1}^{d}r_{i}V^{m}_{1}(V^{*}_{1})^{m}r_{i}^{*} \).
Then, writing
\begin{eqnarray*}
1 & = & V_{1}^{m}(V_{1}^{*})^{m}+(1-V_{1}^{m}(V_{1}^{*})^{m})\\
 & = & V_{1}^{m}(V_{1}^{*})^{m}+\sum r_{i}V_{1}^{m}(V_{1}^{*})^{m}r_{i}^{*}
\end{eqnarray*}
we obtain
\begin{eqnarray*}
W & = & (V_{1}^{*})^{m}\otimes 1\cdot a^{(1)}\otimes a^{(2)}\cdot 1\otimes V_{2}^{k}\\
 & = & (V_{1}^{*})^{m}\otimes 1\cdot a^{(1)}\otimes a^{(2)}\cdot V_{1}^{m}\otimes 1\cdot (V_{1}^{*})^{m}\otimes V_{2}^{k}+\\
 &  & \sum _{i=1}^{d}(V_{1}^{*})^{m}\otimes 1\cdot a^{(1)}r_{i}\otimes a^{(2)}\cdot V_{1}^{m}\otimes 1\cdot (V_{1}^{*})^{m}r_{i}^{*}\otimes V_{2}^{k}\\
 & \in  & \textrm{span}\{(A^{\phi _{1}}_{1}\otimes A_{2}^{\phi _{2}})\cdot (R_{1}\otimes R_{2})^{\phi _{1}\otimes \phi _{2}}\}.
\end{eqnarray*}
Reversing the roles of \( A_{1} \) and \( A_{2} \), we get that in general,
\( \textrm{span}\{(A_{1}^{\phi _{1}}\otimes A_{2}^{\phi _{2}})\cdot (R_{1}\otimes R_{2})^{\phi _{1}\otimes \phi _{2}}\} \)
is dense in \( (A_{1}\otimes A_{2})^{\phi _{1}\otimes \phi _{2}} \).

Since each \( R_{i} \) is hyperfinite, the algebra \( R_{1}\otimes R_{2} \)
is also hyperfinite; hence \( (R_{1}\otimes R_{2})^{\phi _{1}\otimes \phi _{2}} \)
is hyperfinite. It follows that the centralizer \( \phi _{1}\otimes \phi _{2} \)
of \( M=A_{1}\otimes A_{2} \) can be written as the closure of the span of
\( NR \), where \( N \) is a tensor product of two type II\( _{1} \) factors,
and \( R \) is a hyperfinite algebra. Since every hyperfinite algebra can be
written as a linear span of the product \( C_{1}\cdot C_{2} \), where \( C_{i} \)
are abelian von Neumann algebras, it follows that the centralizer \( M^{\phi } \)
is the closure of the span of \( N\cdot C_{1}\cdot C_{2} \), with \( N \)
a tensor product of two type II\( _{1} \) factors, and \( C_{1} \), \( C_{2} \)
abelian von Neumann algebras. Hence by Stefan's result \cite{stephan:thinness},
we get that \( M^{\phi } \) cannot be isomorphic to \( L(\mathbb {F}_{\infty }) \).
\end{proof}
\begin{thm}
Let \( T_{\lambda } \) be the free Araki-Woods factor constructed in \cite{shlyakht:quasifree:big}.
Then \( T_{\lambda }\not \cong A_{1}\otimes A_{2} \), where \( A_{1} \) and
\( A_{2} \) are any diffuse von Neumann algebras.
\end{thm}
\begin{proof}
Since \( T_{\lambda } \) is a full III\( _{\lambda } \) factor, we have by
Lemma \ref{lemma:decomposition}, that the only possible tensor product decompositions
with \( A_{1} \) and \( A_{2} \) diffuse are ones where either exactly one
of \( A_{1} \) and \( A_{2} \) is type II\( _{1} \) and the other is of type
III\( _{\lambda } \), or each \( A_{i} \) is of type III\( _{\lambda _{i}} \),
with \( \lambda _{1}^{\mathbb {Z}}\cap \lambda _{2}^{\mathbb {Z}}=\lambda ^{\mathbb {Z}} \). 

Denote by \( \psi  \) the free quasi-free state on \( T_{\lambda } \). It
is known (see \cite[Corollary 6.8]{shlyakht:quasifree:big}) that \( T_{\lambda }^{\psi } \)
is a factor, isomorphic to \( L(\mathbb {F}_{\infty }) \). Let \( \phi  \)
be an arbitrary normal faithful state on \( T_{\lambda } \), such that \( T_{\lambda }^{\phi } \)
is a factor. Then (see Theorem \ref{thrm:connesresultes}), \( T_{\lambda }^{\phi }\otimes B(\ell ^{2})\cong T_{\lambda }^{\psi }\otimes B(\ell ^{2})\cong L(\mathbb {F}_{\infty })\otimes B(\ell ^{2}) \).
Since \( L(\mathbb {F}_{\infty }) \) has \( \mathbb {R} \) as its fundamental
group (see \cite{radulescu2}), it follows that whenever \( \phi  \) is a state
on \( T_{\lambda } \), and \( T_{\lambda }^{\phi } \) is a factor, then \( T_{\lambda }^{\phi }\cong L(\mathbb {F}_{\infty }) \).

Assume now that one of \( A_{1} \), \( A_{2} \) is of type II\( _{1} \);
for definiteness, assume that it is \( A_{1} \). Choose on \( A_{2} \) a normal
faithful state \( \phi _{2} \) for which \( A^{\phi _{2}} \) is a factor,
and let \( \tau  \) be the unique trace on \( A_{1} \). Let \( \phi =\tau \otimes \phi _{2} \)
on \( T_{\lambda } \). Then \( T_{\lambda }^{\phi }\cong A_{1}\otimes A_{2}^{\phi _{2}} \),
and hence cannot be isomorphic to \( L(\mathbb {F}_{\infty }) \) by a results
of Stephan \cite{stephan:thinness} and Ge \cite{ge:entropy2}, which is a contradiction.

Assume now that \( A_{i} \) is type III\( _{\lambda _{i}} \), with \( 0<\lambda _{i}<1 \).
Then by Proposition \ref{label:prop:centralizer}, there is a state \( \phi  \)
on \( T_{\lambda } \), for which \( T_{\lambda }^{\phi } \) is a factor, but
is not isomorphic to \( L(\mathbb {F}_{\infty }) \), which is a contradiction.
\end{proof}

\section{\protect\( T_{\lambda }\protect \) has no Cartan subalgebras.}

Recall that a von Neumann algebra \( M \) is said to contain a \emph{Cartan
subalgebra \( A \),} if:

\begin{enumerate}
\item \( A\subset M \) is a MASA (maximal abelian subalgebra)
\item There exists a faithful normal conditional expectation \( E:M\to A \)
\item \( M=W^{*}(\mathcal{N}(A)) \), where \( \mathcal{N}(A)=\{u\in M:uAu^{*}=A,\, \, u^{*}u=uu^{*}=1\} \)
is the normalizer of \( A \).
\end{enumerate}
For type II\( _{1} \) factors \( M \), condition (2) is automatically implied
by (1).

\begin{prop}
Let \( M \) be a factor of type \( \III _{\lambda } \), \( 0<\lambda <1 \).
Then there exists a normal faithful state \( \psi  \) on \( M \), so that
\( \sigma ^{\psi }_{2\pi /\log \lambda }=\operatorname {id} \), and that the
centralizer \( M^{\psi } \) is a \( \II _{1} \) factor containing a Cartan
subalgebra.
\end{prop}
\begin{proof}
Let \( A\subset M \) be a Cartan subalgebra. Let \( E:M\to A \) be a normal
faithful conditional expectation. Let \( \phi  \) be a normal faithful state
on \( A\cong L^{\infty }[0,1] \), and denote by \( \theta  \) the state \( \phi \circ E \)
on \( M \). Then \( \theta  \) is a normal faithful state. Furthermore, \( M^{\theta }\supset A \),
because \( E \) is \( \theta  \)-preserving and hence \( \sigma ^{\theta }|_{A}=\sigma ^{\theta |_{A}}=\textrm{id} \).
Since \( M \) is type \( \III _{\lambda } \), it follows that \( \sigma ^{\theta }_{t_{0}} \)
is inner if \( t_{0}=2\pi /\log \lambda  \). Let \( u\in M \) be a unitary
for which \( \sigma ^{\theta }_{t_{0}}(m)=umu^{*} \), \( \forall m\in M \).
Then \( uxu^{*}=x \) for all \( x\in A \), since \( \sigma ^{\theta }|_{A}=\textrm{id} \).
It follows that \( u\in A'\cap M=A' \), since \( A \) is a MASA. Choose \( d\in A \)
positive so that \( d^{it_{0}}=u \). Note that \( d \) is in the centralizer
of \( \theta  \) (which contains \( A \)). Set \( \psi (m)=\theta (d^{-1}m) \)
for all \( m\in M \). Then the modular group of \( \psi  \) at time \( t_{0} \)
is given by \( \textrm{Ad}_{u^{*}}\circ \sigma _{t_{0}}^{\theta }=\textrm{id} \).
It follows that \( \psi  \) is a normal faithful state on \( M \), so that
\( \sigma ^{\psi }_{t_{0}}=\textrm{id} \). It furthermore follows from \ref{thrm:connesresultes}
that the centralizer of \( M^{\psi } \) is a factor of type \( \II _{1} \).
By the choice of \( \psi  \), its modular group fixes \( A \) pointwise, hence
\( A\subset M^{\psi } \).

We claim that \( A \) is a Cartan subalgebra in \( N=M^{\psi } \). First,
\( A'\cap N\subset A'\cap M=A \), hence \( A \) is a MASA. Since \( A \)
is a Cartan subalgebra in \( M \), \( M \) is densely linearly spanned by
elements of the form \( f\cdot u \), where \( u\in \mathcal{N}(A) \) is a
unitary and \( f\in A \). The map
\[
E(m)=\frac{1}{2\pi }\int _{0}^{2\pi }\sigma _{t}^{\psi }(m)dt\]
is a normal faithful conditional expectation from \( M \) onto \( N \). If
\( u\in \mathcal{N}(A) \) is a unitary, so that \( ufu^{*}=\alpha (f) \) for
all \( f\in A \) and \( \alpha \in \textrm{Aut}(A) \), then \( uf=\alpha (f)u \).
Hence
\[
E(u)f=E(uf)=E(\alpha (f)u)=\alpha (f)E(u).\]
It follows that \( N \) is densely linearly spanned by elements of the form
\( E(f\cdot u)=f\cdot E(u) \) for \( f\in A \) and \( u\in \mathcal{N}(A) \).
Let \( w(u) \) be the polar part of \( E(u) \), and let \( p(u)=E(u)^{*}E(u) \)
be the positive part of \( E(u) \), so that \( E(u)=w(u)p(u) \) is the polar
decomposition of \( E(u) \). Since
\[
E(u)^{*}E(u)\alpha ^{-1}(f)=E(u)^{*}fE(u)=\alpha ^{-1}(f)E(u)^{*}E(u),\]
it follows that \( p(u) \) commutes with \( A \) and hence is in \( A \).
Moreover, we then have that
\[
w(u)fw(u)^{*}=\alpha (f),\]
so that \( w(u)\in \mathcal{N}(A)\cap N \) . Thus \( N \) is densely linearly
spanned by elements of the form \( f\cdot u \) for \( f\in A \) and \( u\in \mathcal{N}(A)\cap N \),
hence \( A \) is a Cartan subalgebra of \( N \).
\end{proof}
\begin{cor}
For each \( 0<\lambda <1 \) the \( \III _{\lambda } \) free Araki-Woods factor
\( T_{\lambda } \) does not have a Cartan subalgebra.
\end{cor}
\begin{proof}
If \( T_{\lambda } \) were to contain a Cartan subalgebra, it would follow
that for a certain state \( \psi  \) on \( T_{\lambda } \), the centralizer
of \( \psi  \) is a factor containing a Cartan subalgebra. Let \( \phi  \)
be the free quasi-free state on \( T_{\lambda } \). Then by \ref{thrm:connesresultes},
one has
\[
(T_{\lambda })^{\phi }\otimes B(\ell ^{2})\cong (T_{\lambda })^{\psi }\otimes B(\ell ^{2}).\]
Since \( (T_{\lambda })^{\phi }\cong L(\mathbb {F}_{\infty }) \) (see \cite[Corollary 6.8]{shlyakht:quasifree:big}),
and because the fundamental group of \( L(\mathbb {F}_{\infty }) \) is all
of \( \mathbb {R}_{+} \) (see \cite{radulescu2}) we conclude that \( L(\mathbb {F}_{\infty }) \)
contains a Cartan subalgebra. But this is in contradiction to a result of Voiculescu
that \( L(\mathbb {F}_{\infty }) \) has no Cartan subalgebras (see \cite{dvv:entropy3}).
\end{proof}
\bibliographystyle{amsplain}

\providecommand{\bysame}{\leavevmode\hbox to3em{\hrulefill}\thinspace}

\end{document}